\documentclass{article}

\usepackage{graphicx}
\usepackage{hyperref}
\usepackage{cite}
\usepackage[utf8]{inputenc}
\usepackage{amsmath}
\usepackage{amsfonts}
\usepackage{amssymb}
\usepackage{xcolor}
\usepackage{amsthm}
\usepackage[all]{xy}
\usepackage{enumitem}
\usepackage{tikz-cd}

\usetikzlibrary{cd}

\usepackage{todonotes}
\usepackage{etoc}
\usepackage{amsmath} 

\newtheorem{thm}{Theorem}[section]
\newtheorem{lemma}[thm]{Lemma}
\newtheorem{prop}[thm]{Proposition}
\newtheorem{cor}[thm]{Corollary}

\theoremstyle{definition}
\newtheorem{defn}[thm]{Definition}
\newtheorem{example}[thm]{Example}

\theoremstyle{remark}
\newtheorem{rem}[thm]{Remark}

\numberwithin{equation}{section}

\usepackage{float}

\title{ $\chi$ -extending modular lattices}  
\author{\textbf{Jesus Adrian Celis-González and  Hugo Alberto Rincón-Mejía}}
\date{}

\begin{document}
\maketitle
\begin{multicols}{2}
 \columnbreak

\begin{flushright}\noindent\textbf{MSC2020:}\hspace{1cm}\phantom{.}\\
06C05,06B35,06C15,06C20 \end{flushright}

\end{multicols}

\begin{abstract}
This paper investigates the theory of lattices, focusing on extending lattices relative to abstract classes, modular lattices, and torsion lattices. Definitions of type-1 and type-2 extending lattices are provided, along with their weakly extending counterparts. Modular and upper continuous lattices are analyzed, with lemmas establishing relationships between pseudo-complements, complements, and essential elements. Applications to module theory are highlighted. 
\end{abstract}

  \section{Introduction} A lattice $(L,\leq , \wedge, \vee)$ is bounded if it contains two elements, $0$ (the least element) and $1$ (the greatest element), such that for every element $a \in L$, the inequalities $0 \leq a \leq 1$ are satisfied.
    
 In a lattice \( L \) with the least element \( 0 \), \( b \) is defined as a pseudo-complement of \( a \) if it is a maximal element in the set \( \{ x \in L \mid a \wedge x = 0 \} \). An element \( a \) is considered essential if \( b \neq 0 \) implies \( a \wedge b \neq 0 \). We denote \( b/a = \{ x \in L \mid a \leq x \leq b \} \), calling it an interval in \( L \). In a lattice that includes \( 0 \), an element is (essentially) closed if \( a \) being essential within \( b/0 \) implies \( a = b \). In a bounded lattice, \(b\) is a complement to \(a\) if \(a \vee b = 1\) and \(a \wedge b = 0\). We write $a\oplus b=a\vee b$, when $a\wedge b = 0 $. A lattice with exactly two elements is called simple. In a lattice $\mathcal{L}$, $b$ dominates $a$ if the interval $b/a$ is a simple lattice. A lattice $L$ is complete when each subset $X$ of $L$ possesses a supremum and an infimum. 
\begin{defn}
    A complete lattice $L$ is called upper continuous if  $a\wedge (\bigvee D)=\bigvee_{d\in D} (a\wedge d)$  for every $a\in L$ and for every upper directed subset $D\subseteq L$.
  
\end{defn}
This definition appears in \cite{Calugareanu2000}, Chapter 2, p. 17.
\begin{defn}[Nachbin, Stenstrom]
An element $c$ of a complete lattice $L$ is called \emph{compact} if for every subset $X$ of $L$ and $c \leq \bigvee X$, there is a finite subset $F \subseteq X$ such that $c \leq \bigvee F$, and \emph{$S$-compact} if for each upper directed subset $D \subseteq L$ and $c \leq \bigvee D$, there is an element $d_0 \in D$ such that $c \leq d_0$.
\end{defn}

In Calugareanu's book, it is proven that these two concepts are equivalent.

\begin{defn}
   A lattice $L$ is modular if for each $a, b, c \in L$, the condition $b \leq a$ implies that $a \wedge (b \vee c) = b \vee (a \wedge c)$.
\end{defn}

\begin{defn}
    A complete lattice $L$ is called compactly generated if each of its elements is a join of compact elements.
\end{defn}
\begin{lemma}
    Every compactly generated lattice is upper continuous.
\end{lemma}
\begin{defn}
 An idiom is a complete, upper-continuous modular lattice. 
\end{defn}

\begin{defn}
    Let $L$ be a bounded lattice and $a\in L$. We define the following sets:
$P(a)=\{x\in L:x$ is a pseudocomplement of $a\}$.
     $P(L)=\{x\in L:$ There exists $y\in L$ such that $x\in P(y)\}$.
      $C(L)=\{x\in L: x$ is a closed element in $L\}$.
      $D(L)=\{x\in L:x$ is a complement element in $L\}$.
      $E(L)=\{x\in L:x$ is an essential element in $L\}$.
    $E(L)=\{x\in L:x$ is an essential element in $L\}$.
\end{defn}

    \begin{rem}
    In a bounded modular lattice, $D(L)\subseteq C(L).$

    \end{rem}
\begin{defn}
   A class of lattices $\mathcal{X}$ is abstract if it contains the lattice $0$ and is closed under isomorphic copies of its objects. A $\mathcal X$-lattice $L$, is a lattice belonging to $\mathcal X$. If $\mathcal X$ is a class of lattices, $L$ is a bounded lattice and $a/0$ is an initial interval in $L$, we will say that $a/0$ is a $\mathcal X$-subinterval if $a/0$ is a member of $\mathcal X$.
\end{defn}

\section{Basic concepts}

 \begin{lemma}\label{modular}
    Let $L$ be a modular lattice with least element $0$, and $a,b,c \in L$. If $a\wedge b=0$ and $(a\oplus b)\wedge c=0$, then $a\wedge (b\oplus c)=0$.
\end{lemma}
\begin{proof}
 
   Applying modularity, we can deduce:
\[ a \wedge (b \oplus c) \leq (a \oplus b) \wedge (b \oplus c) = b \oplus ((a \oplus b) \wedge c) = b \oplus 0 = b. \] Since \( a \wedge (b \oplus c) \leq a \), the conclusion follows.
\end{proof}

\
\begin{lemma}\label{caractdeseudocompl}
   Let $L$ be a bounded modular lattice and $a, b \in L$. Then, $b \in P(a)$ if and only if $b \in C(L)$, $a \wedge b = 0$, and $a \oplus b \leq_{e} 1$.
\end{lemma}
\begin{proof}
    
    \(\Longrightarrow)\) Assuming \(b \in P(a)\), we find that \(b \in C(L)\) because condition \(b \leq_e c\) leads to \(a \wedge c = 0\). This follows since \(0 = a \wedge b = b \wedge (a \wedge c)\), and knowing \(b \leq c\) and that \(b \in P(a)\), we deduce \(b = c\). Now, consider \(d \in L\) so that \((a \oplus b) \wedge d = 0\). Referring to Lemma \ref{modular}, it results in \(a \wedge (b \oplus d) = 0\). Given \(b \leq b \oplus d\) and \(b\in P(a)\), we have \(b = b \oplus d\), thus \(d = 0\).

    \(\Longleftarrow)\) To show that $b\in P(a)$, let \(d \in L\) be such that \(b \leq d\) and \(a \wedge d = 0\). We will prove $b=d$, by demonstrating \(b \leq_{e} d\). If \(c \leq d\) with \(c \wedge b = 0\), then \(a \wedge (b \oplus c) \leq a \wedge d = 0\). According to Lemma \ref{modular}, \((a \oplus b) \wedge c = 0\). Given that \(a \oplus b \leq_{e} 1\), we deduce that \(c = 0\). Consequently, \(b \leq_{e} d\), and since \(b \in C(L)\), we conclude \(b = d\).
\end{proof}  
\begin{lemma}\label{caracdeseudocompl2}
    Let $L$ be a bounded modular lattice and $a,b \in L$. $b\in P(a)$ if and only if $a\wedge b=0$ and $a\oplus b \in E(1/b)$.
\end{lemma}
\begin{proof}
  \(\Longrightarrow\)) Let $c\in L$ with $b< c$. Since $c\neq b$, we have $a\wedge c \nleq b$, if $a\wedge c \leq b$, then $a\wedge c \leq a\wedge b=0$, contradicting $b\in P(a)$. Therefore, $b< b\vee (a\wedge c)=(a\vee b)\wedge c$. Thus, $a\oplus b \in E(1/b)$.
  \(\Longleftarrow\) ) Let $c\in L$, be such that $a\wedge c=0$ and $b \leq c$. Suppose $b< c$. Since $a\oplus b \in E(1/b)$, we have $(a\oplus b)\wedge c > b$. On the other hand, using modularity, we have $(a\oplus b)\wedge c=(a\wedge c)\vee b=0\vee b=b$. Therefore, $b<b$, which cannot be. In conclusion, $c=b$.
\end{proof}

\begin{lemma}\label{cerees}
    Let $L$ be a modular and upper continuous lattice. Then, for every $a\in L$, there exists $c\in L $ such that $a\leq_e c $ and $c\in C(L)$.
\end{lemma}
\begin{proof}
   Let $a \in L$. Since $L$ is upper continuous, Zorn's Lemma guarantees the existence of $b \in P(a)$. Similarly, there is $c \in P(b)$ such that $a \leq c$. Because  $c \in P(b)$, $c$  is closed. In fact, if $c \leq_e d$, then $c \wedge (b \wedge d) = b \wedge c = 0$, which implies $b \wedge d = 0$ and, consequently, $c = d$.
Now, we will show that $a \leq_e c$. Let $x \leq c$ be such that $a \wedge x = 0$. Thus, $b\wedge (a\vee x)\leq b \wedge c=0$. By Lemma \ref{modular}, $(b\vee a)\wedge x=0$.  So we have $a\oplus (b\oplus x) $ . Since $b \in P(a)$, according to Lemma \ref{caractdeseudocompl}, then $x=0$. Thus,  $a \leq_e c$.

\end{proof}


\begin{lemma} \label{esenmod1}
  Let $L$ be a lattice with $0$ and $a, b, c \in L$. If $a \leq_e c$, then $a \wedge b \leq_e c \wedge b$. 
\end{lemma} 
\begin{proof} 
If $a \wedge b = 0$, then $b \wedge c = 0$ since $a \leq_{e} c$. Let $x \in c \wedge b / 0$ be nonzero. Since $x \leq c$, we have $a \wedge x \neq 0$ and thus $0 \neq (a \wedge x) = a \wedge (x \wedge b) = x \wedge (a \wedge b)$.
\end{proof}
\begin{lemma}\label{independiente}
    Let $L$ be a modular lattice with a least element, zero, and $a, b, c \in L$. If $a \leq_{e} c$ and $c \wedge b = 0$, then $a \vee b \leq_{e} c \vee b$.
\end{lemma}
\begin{proof}
    Let \( x \leq c \vee b \) be nonzero.   

If \( x \leq b \), then \( 0 < x = x \wedge b \leq x \wedge (a \vee b) \), and hence \( x \wedge (a \vee b) \neq 0 \).  

Now suppose \( x \nleq b \). Then \( x \vee b > b \). Since \( L \) is modular, we have \( (c \vee b)/b \cong c/0 \) via the lattice isomorphism \( c \wedge (-) \). Since \( b < x \vee b \), it follows that \( 0 \neq (x \vee b) \wedge c \). Since \( a \leq_{e} c \), we have \( 0 \neq ((x \vee b) \wedge c) \wedge a = (x \vee b) \wedge a \).  

Assume that \( (a \vee b) \wedge x = 0 \). Since \( a \wedge b = 0 \), by Lemma \ref{modular}, it follows that \( (x \vee b) \wedge a = 0 \), which is a contradiction. Therefore, \( (b \vee a) \wedge x \neq 0 \).  
\end{proof}
 \begin{lemma}
 If $a \wedge b \leq_e b$, then $a \leq_e a \vee b$.
 \end{lemma}
 \begin{proof}Let $a \wedge b \leq_e b$. We want to prove that $a \leq_e a \vee b$. Suppose $f \wedge a = 0$, with $f \leq (a \vee b)$. If $f \neq 0$ then $a < a \oplus f \leq a \vee b$. Applying the lattice isomorphism $b \wedge (-)$ \footnote{Whose inverse is $a\vee(-)$}.
we obtain $b \wedge a < b \wedge (a \oplus f) = (b \wedge a) \oplus (b \wedge f)$. Since $a \wedge b$ is essential in $b$, we obtain $(b \wedge f) = 0$, which contradicts $a \wedge b < b \wedge (a \oplus f)$. Therefore, $f = 0$.\end{proof}
Note that the converse of this result is a consequence of the Lemma \ref{esenmod1} 
\subsection{Linear morphisms of lattices}
The following definition is due to Albu and Iosif in \cite{Albu2013}, chapter 1, page 2.
\begin{defn} \label{linmorf}
    A map between bounded lattices $\varphi : L \rightarrow L'$  is a linear morphism if there exist $k_\varphi\in L$ and $a\in L'$ with the following properties:
is called a linear morphism if there exist $varphi \in L$ and $a\in L'$ such that \begin{enumerate}
   \item \( \varphi(x) = \varphi(x \vee \ker \varphi) \) for all \( x \in L \).
   \item \( \varphi \) induces an isomorphism of the lattices \( \bar{\varphi} : 1/ k_\varphi \to a/0\), given by \( \bar{\varphi}(x) = \varphi(x) \) for all \( x \in 1/ k_\varphi\).
   
\end{enumerate}

\end{defn} 
The following proposition, which can be easily verified, appears in \cite{MBRM}.
\begin{prop} \label{linproj}
Let $L$ be a bounded modular lattice and $x \in L$ be an element with a complement $x'$. The map $\pi_x : L \to L$ is defined by $\pi_x(a) = (a \lor x) \land x'$. It acts as a linear morphism, where $k_{\pi_x} = x'$. The induced isomorphism is $1/x' \to x/0$, and this map is referred to as the $x$ projection. 
\end{prop}
\section{Type-1 and type-2 extending lattices relative to an abstract class of lattices}
\begin{defn}
    A lattice $L$ is extending if for each $a\in L$, there exists $d\in D(L)$ such that $a\in E(d/0)$.
\end{defn}
The following definitions extend the definitions for modules presented in \cite{Dogruöz2}, chapter 1, page 2.
\begin{defn}
    \begin{enumerate}
        Let $L$ be a bounded lattice and let $\mathcal{X}$ be an abstract class of lattices.
        \item A lattice $L$ is type-1 $\mathcal X$-extending if for every $\mathcal X$ interval $a/0$ in $L$, a pseudo-complement $b$ of $a$ in $\mathcal{L}$ belongs to $D(L)$.
        \item $L$ is type-2 $\mathcal{X}$-extending if for every $\mathcal{X}$-interval $a/0$ in $L$, an essential closure of $a$ in $L$ is in $D(L)$.
        \item A lattice $L$ is weakly type-1 $\mathcal X$-extending if for every $\mathcal X$-interval $a/0$ in $L$ there exists a pseudocomplement $b$ of $a$ in $L$, such that $b \in D(L)$.
        \item $L$ is weakly type-2 $\mathcal{X}$-extending if, for every $\mathcal{X}$-interval $a/0$ in $L$, there exists an essential closure of $a$ in $L$ that belongs to $D(L)$.
    \end{enumerate}

\end{defn}

\begin{lemma} \label{seudodeseudo}
    Let $L$ be an idiom. If $a\in L$, $b \in P(a)$, and $c\in P(b)$ are such that $a\leq c, $ then $a \leq_e c$.
\end{lemma}
\begin{proof}
   Since $L$ is an idiom, there exist $b $ and $c $ that satisfy the specified conditions, and both belong to $C(L)$. 
    If $t\leq c$ and $t\wedge a=0$, then $a\wedge (t\vee b)\leq c\wedge (t\vee b)= t\vee (c\wedge b)= t\vee 0=t$.
    
    Thus,  $a\wedge (t\vee b)\leq t\wedge a=0.$ 
    Since $b\in P(a)$ we have $t\vee b=b$. Therefore, $t\leq b$, which implies $t\leq b\wedge c=0$. Hence, $a\leq_e c$.
\end{proof}
\begin{lemma}\label{0.20}
    Let $\mathcal{L}$ be the class of all lattices and $L$ be an idiom. The following statements are equivalent:
\begin{enumerate}
    \item $L$ is extending.
        \item $L$ is type-$1$ $\mathcal{L}$-extending.
        \item $L$ is type-$2$ $\mathcal{L}$-extending.    
\end{enumerate}
\end{lemma}
\begin{proof}

$1) \Longrightarrow 2)$ Let $a/0 \in L$ and $b\in P(a)$. $L$ extending, implies that there exists $d\in D(L)$ such that $b\leq_{e} d$. Since $b\in P(a)$ and $b\leq_{e} d$, by modularity, $a\wedge d=0$. Therefore, $b=d\in D(L)$. 

$2) \Longrightarrow 3)$ Let $a/0 \subseteq L$ and $c\in C(L)$ with $a\leq_{e} c$. Since $L$ is an idiom, there exists $d\in P(a)$. It is clear that $d\wedge c=0$, and there exists $d' \in P(d)$ such that $c\leq d'$. By Lemma \ref{seudodeseudo}, $c\leq_{e} d'$. Thus, $c=d'$ because $c$ is a closed element. So $c=d' \in P(d) \cap D(L)$, since $d' \in P(d)$ and $L$ is type-$1$ $\mathcal{L}$-extending. As $d' \in P(d)$ and $L$ is type-$1$ $\mathcal{L}$-extending, we conclude that $c=d' \in P(d) \cap D(L)$.

$ 3) \Longrightarrow 1) $  Let $ a\in L$, let us prove that there exists $b\in P(a)$ such that $b\in D(L)$. By Lemma \ref{cerees}, there exists $c\in L$ such that $a\leq _e c$ and $c\in C(L)$.  Since $L$ is type-$2$ $\mathcal L$-extending, it follows that $c\in D(L)$. Now, if $c' \in L$ is such that $c\oplus c' = 1$, then $a\wedge c'\leq c\wedge c'=0$. 
If $c' \leq  d$ and $d\wedge a =0$, then $$d=d\wedge 1=d\wedge (c\oplus c')=(d\wedge c)\oplus c',$$ by modularity.
Now, as $a\leq_e c$ and $(d\wedge c)\wedge a=d\wedge a=0$, we obtain $d\wedge c=0$. Hence, $d=c' $. Then $c'\in P(a)\cap D(L)$.

\end{proof}

\begin{lemma} \label{paraidiomastipo1esequivadebtipo1}
   Let $\mathcal{X}$ denote an abstract class of lattices and $L$ be an idiom. If $L$ is a lattice extending type-1 (type-2) $\mathcal{X}$, then it is also a weakly extending type-1 (type-2) $\mathcal{X}$.
 \end{lemma}
 \begin{proof} This is an immediate consequence of Lemma \ref{cerees}.
     
 \end{proof}
\begin{lemma} \label {0.6}
   Let $\mathcal{X}$ denote an abstract class of lattices, and let $L$ be a bounded modular and upper continuous lattice. If $L$ is type-2 $\mathcal{X}$-extending, then $L$ is also weakly type-1 $\mathcal{X}$-extending.
\begin{proof}
Let $a/0$  be an $\mathcal{X}$-interval in $L$. By Lemma \ref{cerees}, there exists $c\in C(L)$ with $a\leq_{e} c$. As $L$ is type-2 $\mathcal{X}$-extending, there exists $c'\in L$ such that $c\oplus c'=1$.  Since $a\leq_{e} c$,  it follows that $a\oplus c' \leq_{e} 1$.  Hence $c' \in P(a)\cap D(L)$.
\end{proof}
\end{lemma}

\begin{lemma}\label{wtype1extpasaasubrets}
Let $\mathcal{X}$ and $\mathcal{Y}$ be abstract classes of lattices such that $\mathcal{X} \subseteq \mathcal{Y}$. If $L$ is a weakly type-$1$  $\mathcal{Y}$-extending,  then  $L$ is weakly type-$1$ $\mathcal{X}$-extending.
\end{lemma}
\begin{proof}
    Let $a\in L$ such that $a/0$ is in $\mathcal{X}$. Since $\mathcal{X} \subseteq \mathcal{Y}$, $a/0$ is $\mathcal{Y}$. Therefore, there exists $b\in P(a)$ such that $b\in D(L)$.
\end{proof}
\begin{defn}
Let $\mathcal{X}$ be an abstract class of lattices having a  least element. We define $\mathcal{X}^{e}$ as the class of lattices $L$ that contain an  $\mathcal{X}$-interval $a/0$ for some $a \in E(L)$. 
\end{defn}
Note that $\mathcal{X}\subseteq \mathcal{X}^e$ if $\mathcal{X}$ consists of bounded lattices.
\begin{prop}
   Let $\mathcal{X}$ be a class of bounded lattices and let $L$ be a modular and bounded lattice. Then $L$ is weakly type-1 $\mathcal{X}$-extending if and only if  $L$ is weakly type-1 $\mathcal{X}^{e}$-extending $1$.
   \end{prop}
\begin{proof}
   $\Longrightarrow )$ Let  $a/0$ be an $\mathcal{X}^{e}$-interval in $L$. Since  $a/0 \in \mathcal{X}^{e}$, there exists an  $\mathcal{X}$-interval $b/0$ in $L$ such that  $b\leq_{e} a$.
    \\
    Since $L$ is weakly type-1  $\mathcal{X}$-extending, there exists $c\in  D(L)$  such that $c\in P(b)$. Since $b\wedge c=0$ and $b\leq_{e} a$, we have $a\wedge c=0$. Moreover, by Lemma \ref{caractdeseudocompl}, $b\oplus c\in E(L)$, which implies that $a\oplus c \in E(L)$. Applying Lemma \ref{caractdeseudocompl} again, we conclude that $c\in D(L)\cap P(a)$. Thus,  $L$ is weakly type-1  $\mathcal{X}^{e}$-extending.
    
   $\Longleftarrow )$ This follows  from Lemma   \ref{wtype1extpasaasubrets}
\end{proof}

\begin{lemma} \label{0.8}
  Let $\mathcal{X}$ be a class of lattices and let $L$ be an idiom; then $L$ is type-$2$ $\mathcal{X}^{e}$-extending if and only if $L$ is weakly type-$2$ $\mathcal{X}^{e}$-extending.
\begin{proof}
  $\Longrightarrow )$ This follows from Lemma \ref{paraidiomastipo1esequivadebtipo1}. 
    \\
    $\Longleftarrow )$ Conversely, assume $L$ is weakly type-$2$ $\mathcal{X}^{e}$-extending. Consider $a/0 \in \mathcal{X}^{e}$, with $a\in L$, and let $b$ be an essential closure of $a$ in $L$. \\ Since $a/0$ is in $\mathcal{X}^{e}$, there exists $c/0$ in $\mathcal{X}$ such that $c\leq_{e} a$. Consequently, $c\leq_{e} b$, so $b/0$ is in $\mathcal{X}^{e}$. Thus, there exists $d\in D(L)\cap C(L)$ such that $b\leq_{e} d$. But $b\in C(L)$ so $b=d$ and then $b\in D(L)$. We conclude that $L$ is type-$2$ $\mathcal{X}^{e}$-extending.

\end{proof}
\end{lemma}
\begin{lemma}
       Let $\mathcal{X}$ be a class of bounded lattices and $L$ be a bounded modular lattice. $L$ is  type-$1$  $\mathcal{X}^{e}$-extending if and only if $L$ is  type-$1$  $\mathcal{X}$-extending. 
\begin{proof}
   $\Rightarrow)$  Let $a/0$ be an  $\mathcal{X}^{e}$-interval in $L$ and $b \in P(a)$. Since $a/0$ is in $\mathcal{X}^{e}$, there exists $c/0$, an $\mathcal{X}$-interval in $L$ with $c \leq_{e} a$.  Since $b \in E(L)$ and $c\wedge b=0$,  it follows that $b \in P(c)$ by Lemma \ref{caractdeseudocompl} .
    Since $L $ is type-$1$ $\mathcal{X}$-extending, we conclude that $b \in D(L)$.
      $\Leftarrow) $ This follows from Lemma \ref{paraidiomastipo1esequivadebtipo1}.
uni\end{proof}
\end{lemma}
\begin{thm}
    Let $\mathcal{X}$  be an abstract class of lattices and $L$ be a modular and upper continuous lattice. Then $L$ is type-2 $\mathcal{X}$-extending if and only if   $L$ is weakly type-2  $\mathcal{X}^{e}$-extending.
\end{thm}
\begin{proof}
    It follows from Lemma \ref{paraidiomastipo1esequivadebtipo1} and Lemma \ref{wtype1extpasaasubrets}. 
\end{proof}

We refer the reader to \cite{Calugareanu2000}  for the definition and properties of Goldie finite dimension lattices.

$\mathcal{L}$: the class of all lattices.\\$\mathcal{C}$: the class of all simple lattices.\\$\mathcal{U}$: the class of all lattices with finite Goldie dimension.\\$\mathcal{U}_{1}$: the class of uniform lattices.\\$\mathcal{G}$: the class of compact lattices.

\begin{defn}

    $L$ is extending if for every $a\in L$ there exists $d\in D(L)$ such that  $a\in E(d/0)$.
\end{defn}

\begin{example} In the example below, the lattice $L$ is (weakly)  type-$1$ ($2$) $\mathcal{G}$-extending, but it is not an extending lattice. Here, the compact elements are the natural numbers. It's important to note that $i$ is not a compact element since $i\leq \omega=\vee\{n\in \mathbb{N}\}$, yet $i\nleq n$ for any $n\in \mathbb{N}$. A $\mathcal{G}$-initial interval within $L$ takes the form $n/0$, and the only pseudocomplement of $n$ is $i$, with  $i\in D(L)$. Thus, ${L}$ is indeed type-$1$ $\mathcal{G}$-extending. However, $L$ is not extending because $i$ lacks a pseudocomplement in $D(L)$.
   \begin{figure}[H]
  \centering
  \includegraphics[scale=0.25]{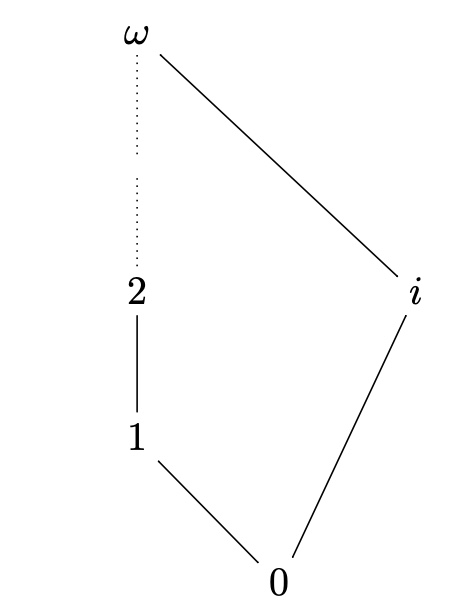}
  \caption{}
  \label{fig:enter-label-2}
\end{figure}
\end{example}

\begin{prop}
    A modular lattice $L$  is weakly type-2 $\mathcal{L}$ extending if and only if it is extending.
\end{prop}\begin{proof}
$\Longrightarrow )$ Consider an element $a \in L$. Given that $a/0$ forms an $L$-interval, there exists an element $d \in D(L)\cap C(L)$ such that $a\leq_{e} d$.
\\
$\Longleftarrow )$ If $a \in L$, then there exists $d \in D(L)$ such that $a \leq_e d$. When $L$ is a modular lattice, we have $D(L) \subseteq C(L)$. Note that $d$ is an essential closure of $a$.
\end{proof}

\begin{lemma}  
Let \( L \) be a bounded modular lattice with \( 1 = s \oplus u \), where \( s \in L \) is an atom and \( u \in L \) is uniform. Then \( L \) is weakly type-1 \(\mathcal{L}\)-extending.  
\end{lemma}  
\begin{proof}  
Let \( a \in L \). We show that \( a \) is essential in an element of \( D(L) \), considering two cases:  

\begin{itemize}  
\item If \( a \wedge s \neq 0 \): Since \( s \) is an atom, \( a \wedge s = s \), we have \( s \leq a \).  
    \begin{itemize}  
    \item If \( a \wedge u = 0 \), modularity gives \( a = a \wedge 1 = a \wedge (u \vee s) = (a \wedge u) \vee s = s \). Thus, \( u \in P(a) \cap D(L) \).  
    \item If \( a \wedge u \neq 0 \) and \( u \) is uniform, then \( a \wedge u \leq_e u \). By modularity, \( \frac{(a \vee u)}{a} \) is isomorphic to \( \frac{u}{(a \wedge u)} \), which implies \( a \leq_e (a \vee u) \). Since \( 1 = s \oplus u \leq (a \vee u) \), it follows that \( a \leq_e 1 \), and thus \( 0 \in P(a) \cap D(L) \).   
    \end{itemize}  

\item If \( a \wedge s = 0 \):   
    \begin{itemize}  
    \item If \( a \neq 0 \), then \( (a \oplus s) \wedge u = 0 \) implies, by Lemma \ref{modular}, \( a = a \wedge (s \oplus u) = 0 \), a contradiction. Thus, \( (a \oplus s) \wedge u \neq 0 \), so \( (a \oplus s) \wedge u \leq_e u \). By Lemma \ref{independiente}, \( (a \oplus s) \leq_e (a \oplus s) \vee u = 1 \). By Lemma \ref{caractdeseudocompl}, \( s \in P(a) \cap D(L) \).   
    \item If \( a = 0 \), then \( 1 \in D(L) \cap P(a) \).  
    \end{itemize}  
\end{itemize}  
\end{proof}

\begin{lemma}\label{cerrados}
   Let $L$ be a bounded modular, $a \in C(L)$, and $b \in C(a/0)$. Then $b \in C(L)$.
\end{lemma}
\begin{proof}

Let \( d \in L \) be such that \( b \leq_{e} d \).
\begin{itemize}
    \item 

If \( a \leq d \), then \( b \leq a \leq d \) with \( b \leq_{e} d \), which implies \( a \leq_{e} d \). Thus, \( a = d \), and consequently, \( b \leq_{e} a \) and \( b \in C(a) \), leading to \(d= a = b \).

\item Suppose \( a \nleq d \). Then \( b \leq a \wedge d < a < a \vee d \), and since \( b \leq_{e} d \), it follows that \( b \leq_{e} (a \wedge d) \). Given that \( b \in C(a/0) \), we conclude \( a \wedge d = b \). Thus, \( (a \wedge d) \leq_{e} d \). By modularity, \( a \leq_{e} (a \vee d) \). This implies \( a = a \vee d \) since \( a \in C(L) \). Therefore, \( b\leq_e d \leq a \), and as \( b \in C(a/0) \), this leads to \( b = d \).
\end{itemize}

\end{proof}
\begin{lemma}\label{pseudorela}
    Let $L$ be an  idiom such that $1=a\oplus b$ and $c,d \leq a$. $d$ is a relative pseudocomplement of  $c$ in $a/0$ if and only if $d\oplus b \in P(c)$.
\end{lemma}
\begin{proof}

 $\Longrightarrow )$ Since $d$ is a relative pseudocomplement of $c$ in $a/0$. It follows that $d\oplus c\leq_{e} a$, and  $(d\oplus c)\oplus b\leq_{e} 1$. By Lemma \ref{modular},  $c\oplus (d\oplus b)\leq_{e} 1$. Therefore, $d\oplus b\in P(c)$, by Lemma \ref{caractdeseudocompl}.
 
$\Longleftarrow )$ Since $c\oplus (d\oplus b)\leq_{e} 1$, then by Lemma \ref{esenmod1}, $(c\oplus (d\oplus b))\wedge a \leq_{e} a$. Now, $(c\oplus (d\oplus b))\wedge a=c\oplus (a\wedge (d\oplus b))=c\oplus d$. Thus, $d$ is a relative pseudocomplement of $c$ in $a/0$.
\end{proof}

\begin{lemma}
   Let $L$ be an idiom and $\mathcal{X}$ a class of lattices. $L$ is type-$1$ $\mathcal{X}$-extending if and only if, for all $a \in L$, the quotient $a/0$ is type-$1$ $\mathcal{X}$-extending.
\end{lemma}
\begin{proof}
$\Longrightarrow )$ Let $a\in D(L)$ and $b\leq a$ be such that  $b/0 \in \mathcal{X}$.  Let $c$ be a pseudocomplement of $b$ in $a/0$. If $a' \in L$ is such that  $a\oplus a'=1$, then $c\oplus a' \in P(b)$,  by Lemma \ref{pseudorela}. 
\\
Since $L$ is type-$1$ $\mathcal{X}$- extending, we have $c\oplus a'\in D(L)$. Let $d' \in L$ be such that $(c\oplus a') \oplus d'=1$, then $a= a\wedge 1=a \wedge ((c\oplus a')\oplus d')=a \wedge (c\oplus (a' \oplus d'))=c \oplus (a\wedge (a'\oplus d'))$,  by modularity. Therefore, $c\in D(L)$.\\
$\Longleftarrow )$ It results from $L=1/0$.
\end{proof}

\begin{lemma}
    Let $L$ be an idiom and  $\mathcal{X}$ be a class of lattices. $L$ is type-$2$ $\mathcal{X}$-extending if and only if for all $c\in C(L)$, $c/0$ is type-$2$ $\mathcal{X}$-extending.
\end{lemma}
\begin{proof}
   $\Longrightarrow )$ Let $c \in C(L)$ and $a \leq c$ such that $a/0 \in \mathcal{X}$. Suppose $d \in C(c/0)$ satisfies $a \leq_{e} d$. By Lemma \ref{cerrados}, we have $d \in C(L)$. Given that $L$ is type-$2$ $\mathcal{X}$-extending, it follows that $d \in D(L)$. Hence, $c/0$ is type-$2$ $\mathcal{X}$-extending.
\\
$\Longleftarrow )$ This is deduced from the fact that $1 \in C(L)$ and $L = 1/0$.
\end{proof}

\begin{lemma}\label{ora}
    Let $L$ be a bounded modular lattice and $c\in C(L)$. Let $a>c$ with $a\in E(L)$, then $a\in E(1/c)$.
\end{lemma}
\begin{proof}
    Consider $b \in 1/c$ such that $b > c$. Since $c \in C(L)$, it follows that $c \not\in E(b/0)$. Consequently, there exists a non-zero $d$ with $d \leq b$ and $c \wedge d = 0$. Observe that $a \wedge (c \vee d) = c \vee (a \wedge d) > c$; otherwise, $a \wedge d \leq c \wedge d = 0$, which would contradict $a \in E(L)$. Thus, $a \wedge b \geq a \wedge (c \vee d) > c$, implying that $a \in E(1/c)$.
    \end{proof}
    
   In the proof of the following theorem, we use the following straightforward observation: If $b \in P(a)$ in a bounded modular lattice $L$, then $a \vee (\cdot):b/0 \rightarrow (a \vee b)/a$ is a lattice isomorphism whose inverse is $b \wedge (\cdot):(a \vee b)/a \rightarrow b/0$.

\begin{thm}
  Let \(\mathcal{X}\) be a class of lattices closed under initial intervals. An idiom \(L\) is type-1 \(\mathcal{X}\) extending if and only if for any \(a \in C(L)\), the existence of \(c \in E(1/a)\) such that \(c/a \in \mathcal{X}\) implies \(a \in D(L)\).

\end{thm}
\begin{proof}
    $\Longrightarrow )$ Let   $a\in C(L)$. Assume that there exists $c\in E(1/a)$ with $c/a  \in \mathcal{X}$. Suppose $b\in P(a)$. Given that $b \oplus a \leq_{e} 1$ and $a \in  C(L)$, from Lemma \ref{caractdeseudocompl} it follows that $a$ is also in $P(b)$.
      On the other hand,  consider the lattice morphism $\pi: L \longrightarrow 1/a$, defined by  $\pi(x)=x\vee a$. Since $b\in P(a)$, it follows that $\pi|_{b/0}$ is an injective lattice morphism. Therefore, $\pi(b)/a = \pi [b/0] \cong b/0$. Given that $c\in E(1/a)$ , it follows that $c \in E(\pi(b)/a)$. Moreover, $(\pi(b)\wedge c )/a\in \mathcal{X}$, since $\mathcal{X}$ is closed under taking initial intervals. In conclusion, $b/0$ has an essential element $e$ such that $e/0\in \mathcal{X}$.

Since $a\in P(b)$, it follows that $a\in P(e)$. Since $L$ is type-1 $\mathcal{X}$-extending, $a\in D(L)$.

$\Longleftarrow )$ Let $a\in L$ such that $a/0 \in \mathcal{X}$ and $c\in P(a)$. By modularity, we have 
\begin{center}
    $a/0=a/(a\wedge c)\cong (c\oplus a)/c$.
\end{center}
Therefore, $(c\oplus a)/c\in \mathcal{X}$. Furthermore, by Lemma \ref{ora}, $c\oplus a\in E(1/c)$. We conclude that $(c\oplus a)\in D(L)$. Thus $c\in D(L)$.
\end{proof}
\begin{prop}

    Let $L$ be an idiom. $L$ is of type-$1$ $\mathcal{U}$-extending if and only if for all $a\in C(L)$ such that $1/a$ has finite uniform dimension, it follows that $a\in D(L)$.

\end{prop}
\begin{proof}
   $\Longrightarrow )$ Suppose $a\in C(L)$,  $a/0 \in \mathcal{U}$, and $b\in P(a)$. Then $b\in C(L)$ and $(a\oplus b) \leq_e 1$. As $(a\oplus b)/b \cong a/0$, and $a/0 \in \mathcal{U}$, then $(a\oplus b)/b \in \mathcal{U}$. Also, $a\oplus b\in E/(L)$ implies $1/b \in \mathcal{U}$. Hence, by the actual hypothesis, $b\in D(L)$

   $\Longleftarrow )$ Suppose now that $L$ is of type$1$ $\mathcal{U}$-extending and $1/a \in \mathcal{U}, a\in C(L)$. Take $b\in P(a)$, which exists because $L$ is an idiom. Then $a\oplus b \leq_e 1$, and since $a\in C(L)$ then $a \in P(b)$.  $b/0\cong (a\oplus b)/a \in \mathcal{U}$ implies that $1/a \in \mathcal{U}$ (because if $\{d_i\}_I$ were an infinite independent set of non-zero elements of $L$ then $\{(a\oplus b)\wedge d_i\}_I$ would also be an infinite independent set.  As $b/0 \in \mathcal{U}$ with $b\in C(L)$, then $a\in D(L)$.
\end{proof}
Recall that a compactly generated lattice is upper continuous, (see \cite{Calugareanu2000}, Theorem  2.4). Thus each  element in a compactly generated lattice has a Pseudocomplement.
\begin{prop}
    Let $L$ be a compactly generated, modular, and indecomposable lattice. Then the following statements are equivalent:
    \begin{enumerate}
        \item $L$ is uniform.
        \item $L$ is type-$1$ $\mathcal{G}$-extending.
        \item $L$ is type-$2$ $\mathcal{G}$-extending.
    \end{enumerate}
\end{prop}
\begin{proof}
    $1)\Longrightarrow 2)$
Let $x\in L$. Since $L$ is uniform, for $x\neq 0$, the sole pseudocomplement of $x$ is $0$,  and $0\in D(L)$, while for $x=0$, its unique  pseudocomplement is $1$, and $1\in D(L)$.
     \\
  $2)\Longrightarrow 3)$ 
Let $x\in L\setminus\{0\}$ be compact and $c\in C(L)$ such that $x\leq_{e} c$. Let $d\in L$ be a pseudocomplement of $c$ and $t\in L$ be a pseudocomplement of $d$ such that $c\leq t$. Let us show that $c\leq_{e} t$. Let $y\leq t$ such that $c\wedge y=0$. By Lemma \ref{modular}, $c\wedge (d\vee y)=0$. Since $d\in P(c)$, we have that $d\vee y=d$ so $y=0$. Since $c\leq_{e} t$ and $c\in C(L)$, we conclude that $c=t$.
    \\
    \(3)\Longrightarrow 1)\) Let   \( a \in L \setminus \{0\} \). Since \( L \) is compactly generated, there exists a compact element \( x \in L \setminus \{0\} \) such that \( x \leq a \). Consider \( c \in C(L) \) such that \( x \leq_{e} c \). Because \( L \) is indecomposable and type-\(2\) \(\mathcal{G}\)-extending, it follows that \( c = 1 \). Consequently, we have \( x \leq_{e} 1 \), which, in particular, implies that \( a \leq_{e} 1 \).
    \end{proof}
\section{\texorpdfstring{$\displaystyle\bigoplus_{i=1}^{n}\mathcal{X}_{i}$-extending lattices}{Direct sum Xi-extending lattices}}    

\begin{defn}

Let \(\mathcal{X}_i\) denote a family of lattice classes for \(1 \leq i \leq n\). The expression \(\mathcal{X}_1 \oplus \ldots \oplus \mathcal{X}_n\) represents the class of lattices \(L\) such that there exists an independent set \(\{a_1, \ldots, a_n\} \subseteq L\) with \(L = \bigoplus_{i=1}^{n} a_i / 0\), where each \(a_i / 0\) belongs to the class \(\mathcal{X}_i\) for all \(1 \leq i \leq n \)

\end{defn}

   \begin{thm}\label{sumadirectatipo1}
   Let $L$ be an idiom and $\mathcal{X}$ be a collection of lattices. $L$ is  type-$1$  $\displaystyle\bigoplus_{i=1}^{n}\mathcal{X}_{i}$-extending if and only if $L$ is  type-$1$  $\mathcal{X}_{i}$-extending for every $1 \leq i \leq n$.
\end{thm}
\begin{proof}

   $\Longrightarrow )$ This follows directly from Lemma \ref{wtype1extpasaasubrets}. 
   
$\Longleftarrow )$ It suffices to prove it for $n=2$. Suppose $L$ is of type-$1$ $\mathcal{X}_{i}$-extending for $i=1, 2$. Let $a_{1}, a_{2}\in L$ be such that $a_{i}/0$ is in $\mathcal{X}_{i}$ for $i=1, 2$ and $a_{1}\wedge a_{2}=0$. Let $a=a_{1}\oplus a_{2}$ and $b\in P(a)$. We claim $b\in D(L)$.
\\
Since $a_{1}\wedge a_{2}=0$,  Lemma \ref{modular} gives $a_{1}\wedge (a_{2}\oplus b)=0$. As $L$ is an idiom, there exists $c\in P(a_{1})$ such that $a_{2}\oplus b \leq c$. Since $L$ is type-$1$ $\mathcal{X}_{1}$-extending, there exists $c'\in L$ such that $c\oplus c'=1$. By modularity,  $(b\oplus a_{2} \oplus a_{1})\wedge c=(c\wedge a_{1})\oplus (b\oplus a_{2})=b\oplus a_{2}$. If  $t\leq c$ and  $(b\oplus a_{2})\wedge t=0$, then $(((b\oplus a_{2})\wedge c)\wedge t=0$. By Lemma \ref{modular},  $(b\oplus t)\wedge (a_{2}\oplus a_{1})=0$, then, and since $b\in P(a)$,  we get $b\oplus t=b$, so $t=0$. Therefore, $b\oplus a_{2} \leq_{e} c$.

Since $b\in P(a)$, we have $b\in C(L)$,  it follows that $b\in C(c/0)$ and $b\oplus a_2 \leq_{e} c$; thus, $b$ is a pseudocomplement of $a_{2}$ with respect to $c/0$. By Lemma \ref{pseudorela}, we have $b\oplus c' \in P(a_{2})$. Moreover, since $L$ is type-1 $\mathcal{X}_{2}$ extending, there exists $b'\in L$, such that $(b\oplus c')\oplus b'=1$. Thus, $b\in D(L)$,  and so $L$ is type-$1$ $\mathcal{X}_{1}\oplus \mathcal{X}_{2}$-extending.

\end{proof}
  \begin{thm}\label{sumadirectatipo2} 
   Let $L$ be an idiom and $\mathcal{X}$ be a collection of lattices. $L$ is  type-$2$ $\displaystyle\bigoplus_{i=1}^{n}\mathcal{X}_{i}$-extending if and only if $L$ is   type-$2$ $\mathcal{X}_{i}$-extending for every $1 \leq i \leq n$.
   \end{thm}
   \begin{proof}
   $\Longrightarrow$)  This follows directly from Lemma \ref{wtype1extpasaasubrets}.
   
   $\Longleftarrow$)  Now, suppose $L$ is type-$2$ $\mathcal{X}_{i}$-extending for $i=1, 2$.  It suffices to consider the case $n=2$.  Let $a_{i}/0$ in $\mathcal{X}_{i}$ for $i=1,2$ be such that $a_{1}\wedge a_{2}=0$. Let $a=a_{1}\oplus a_{2}$ and $c\in C(L)$ be such that $a \leq_{e} c$.

Since $L$ is an idiom, there exists $b\in C(c/0)$ such that $a_{1}\leq_{e} b$. By Lemma \ref{cerrados}, we have $b\in C(L)$. Given that $L$ is type-$2$ $\mathcal{X}_{1}$-extending, there exists an element $b' \in L$ such that $b' \oplus b = 1$. Consequently, we have $c = c \wedge 1 = c \wedge (b' \oplus b) = b \oplus (c \wedge b')$. Consider the linear morphism $\pi: c/0 \rightarrow (c \wedge b')/0$  defined by $\pi(t) = (t\vee b) \wedge (c \wedge b')$ (see Proposition \ref{linproj}). 

We claim $\pi(a_2) \leq_e c \wedge b'$. Suppose $0 \neq x \leq c \wedge b'$. Since $a_1 \oplus a_2 \leq_e c$, it follows that $0 \neq x \wedge (a_1 \oplus a_2) \leq x \wedge (a_2 \vee b)$. 

Therefore, we can write:
\[
\pi(a_2) \wedge x = ((a_2 \vee b) \wedge (c \wedge b')) \wedge x = (a_2 \vee b) \wedge x \neq 0.
\]
Thus, we conclude that $\pi(a_2) \leq_e c \wedge b'$.

  According to the definition of a linear morphism (see Definition \ref{linmorf}), we have $k_\pi = b$. Consider restricting $\pi$ to $a_2/0$. Given that $a_1 \leq_e b$, it follows that $a_2 \wedge b = 0$. If $0 = \pi(x) = \pi(x \vee b)$, this implies that $x \vee b = b$, leading to $x \leq b$. Therefore, if it follows that $x \leq a_2$, then $x = 0$. Consequently, the restricted linear morphism $\pi_{\mid a_2/0}: a_2/0 \rightarrow \pi(a_2)/0$ satisfies $k_{\pi_{\mid a_2/0}} = 0$, making it a lattice isomorphism.
  
Since $a_{2}/0\cong \pi(a_{2})/0$, it follows that $\pi(a_{2})/0$ belongs to $\mathcal{X}_{2}$. Given  $\pi(a_{2})\leq_{e} c\wedge b'$, and  $c\wedge b' \in C(L)$, along with $L$ being type-2 $\mathcal{X}_{2}$ extending, we conclude that $c\wedge b'\in D(L)$. By modularity, $c \wedge b'$ is in $D(b'/0)$, implying that there exists some element $t \in L$ such that $b' = (c \wedge b') \oplus t$. Consequently, $c \oplus t = 1$, showing that $L$ is a type-2 extending lattice.      
   \end{proof}
   
\begin{defn}
    Let $\mathcal{X}$ be a class of lattices. $\mathcal{X}^{\oplus}$ denotes the class of lattices $L$ such that there exist $a_{1},...,a_{n}\subseteq L$ independent such that $a_{i}/0 \in \mathcal{X}$ for all $1\leq i \leq n$ and $L= \displaystyle\oplus_{i=1}^{n}a_{i}/0$.
\end{defn}
\begin{cor}\label{sumadirectaigual}
    Let $L$ be an idiom and $\mathcal{X}$ be a class of lattices. Then $L$ is type-$1$ (respectively, type-$2$) $\mathcal{X}$- extending if and only if $L$ is type-$1$ (type-$2$) $\mathcal{X}^{\oplus}$-extending.
\end{cor}
\begin{proof}
    This follows from Theorems \ref{sumadirectatipo1} and \ref{sumadirectatipo2}.
\end{proof}
\begin{prop}
   Let $L$ be an idiom. Then $L$ is type-$1$ (type-$2$) $\mathcal{U}$-extending if and only if $L$ is type-$1$ (type-$2$) $\mathcal{U}_{1}$-extending.
\end{prop}
\begin{proof}
    $\Longrightarrow )$ It follows from Lemma \ref{wtype1extpasaasubrets}, since $\mathcal{U}_{1} \subseteq \mathcal{U}$.

$\Longleftarrow )$ Suppose $L$ is a type-$1$ (type-$2$) $\mathcal{U}_{1}$-extending . According to Corollary \ref{sumadirectaigual}, $L$ is also a type-$1$ (type-$2$) $\mathcal{U}_{1}^{\oplus}$-extending . Lemma \ref{0.8}  further implies that $L$ is a type-$1$ (type-$2$) $(\mathcal{U}_{1}^{\oplus})^{e}$-extending. Since $\mathcal{U}$ is a subset of $(\mathcal{U}_{1}^{\oplus})^{e}$, Lemma \ref{wtype1extpasaasubrets} ensures that $L$ is a type-$1$ (type-$2$) $\mathcal{U}$-extending .
\end{proof}
\begin{thm}
   Let $L$ be a bounded modular lattice. Suppose that $L$ and that every direct summand of $L$ is weakly type-1 $\mathcal{U}_{1}$-extending. Then $L$ is weakly type-1 $\mathcal{U}$-extending, and the same holds for all its direct summands.
\end{thm}
\begin{proof}
    Let $a\in L$ be such that $a/0\in U$; that is, $u-dim(a/0)=n$, let us proceed by induction on $n$.
 If $u-dim(a/0)=1$, then $a\in L$ is uniform. Since  $L$  is weakly type-1  $\mathcal{U}_{1}$- extending, we have that there exists $d\in D(L)\cap P(a)$.
     Suppose $u-dim(a/0)=n > 1$. Let $u \leq a$ be uniform, since $L$ is weakly type-1 $U_{1}$ extending, there exist $b,b'\in L$ such that   $1=b\oplus b'$ and   $b \oplus u \leq_{e} 1$.  Let us note that $0=b\wedge u = (b\wedge a)\wedge u$, so $u-dim(b\wedge a/0) < u-dim(a/0)$. By the induction hypothesis, applied to $(b\wedge a)/0$ in the interval $b/0$, we see that there exist $c,c' \in b/0$ such that $b=c\oplus c'$ and $b\wedge (c\oplus a)=c\oplus (b\wedge a) \leq_{e} b$.

 Given that $b\oplus u \leq_{e} 1$, it follows that $(b\oplus u)\wedge b'\leq_{e} b'$. Consider the linear morphism $\pi: L \longrightarrow b'/0$  given by   $\pi(x)=(x\vee b)\wedge b'$ (see Proposition \ref{linproj}). Now, $\pi|_{u/0}$ is an injective linear morphism , since $k_{\pi|_{u/0}}\leq u \wedge b=0$. Therefore,
\begin{center}
$u/0\cong \pi(u)/0=(u\vee b)\wedge b'/0$.
\end{center}
Thus, $ (u\vee b)\wedge b'\in L$ is a uniform element and $(u\vee b)\wedge b' \leq_{e} b'$, then $b'/0$ is a uniform interval.
\\
Suppose $ (a \oplus c)\wedge b' \neq 0$ then $(c\oplus a) \wedge b' \leq_{e} b'$. Since $( b\wedge (c\oplus a)\leq_{e} b$, by modularity, we have $(b\wedge (c\oplus a)) \vee ((c\oplus a)\wedge b') \leq_{e} b\oplus b'=1$. Now, $(b\wedge (c\oplus a)) \vee ((c\oplus a)\wedge b') = (c\oplus a)\wedge (((c\oplus a)\wedge b)\vee b')\leq (c\oplus a)\wedge (b\oplus b')=c\oplus a$. In conclusion, $a\oplus c\leq_{e} 1$ and $c\in D(L)$ because $c\oplus (c'\oplus b')=1$. Thus, $c\in D(L)\cap P(a)$.

On the other hand, if $(a \oplus c)\wedge b'=0$ we have $((c\oplus a)\wedge b)\vee b' \leq_{e} b\oplus b'=1$. Observe that $((c\oplus a)\wedge b)\vee b' \leq (c\oplus a)\vee b'$, so $a\oplus (c\oplus b')\leq_{e} 1$. Given that $(c\oplus b')\oplus c'=1$, we conclude that $b' \oplus c\in D(L)\cap P(a)$.
\end{proof}
\section[(X1 ... Xn)-extending lattices]{$(\mathcal X_1....\mathcal X_n)$-extending lattices}

\begin{defn}

Let \( n \) be a positive integer and let \( \mathcal{X}_{i} \) be a class of lattices, \( 1 \leq i \leq n \). We define the class \( \mathcal{X}_{1} \mathcal{X}_{2} \ldots \mathcal{X}_{n} \) as the class of lattices \( L \) such that there exists a chain of elements
\[
0 = a_{0} \leq a_{1} \leq \ldots \leq a_{n} = 1
\]
such that \( \frac{a_{i}}{a_{i-1}} \in \mathcal{X}_{i} \) for \( 1 \leq i \leq n \). In particular, if \( \mathcal{X}_{i} = \mathcal{X} \) for all \( 1 \leq i \leq n \), we denote the class \( \mathcal{X}_{1} \ldots \mathcal{X}_{n} \) by \( \mathcal{X}^{n} \).

\end{defn}
\begin{thm}
    Let $n$ be a positive integer and $\mathcal{X}_i   (1 \leq i \leq n)$ be classes of lattices such that $\mathcal{X}_i $ is closed under initial intervals or quotients for $2\leq i \leq n$. An idiom $L$ is type-$2$ $\mathcal{X}_{1}...\mathcal{X}_{n}$-extending if and only if $L$ is type-$2$ $\mathcal{X}_{i}$-extending for all $1\leq i \leq n$.
\end{thm}
\begin{proof} 

    $\Longrightarrow )$ Let $1 \leq i \leq n$ and $a/0 \in \mathcal{X}_{i}$. Then $0 \leq 0\leq ...\leq 0 \leq a\leq a \leq ...\leq a$  (with the first $a$ placed at the $i$-eth position) is a chain such that $a/a \cong 0 \in \mathcal{X}_{j}$ for each $i+1 \leq j \leq n$, $a/0 \in \mathcal{X}_{i}$, and $0 = 0/0 \in X_{j}$ for all $0 \leq j \leq i-1$. Therefore, $\mathcal{X}_{i} \subseteq \mathcal{X}_{1}\ldots\mathcal{X}_{n}$, and by Lemma \ref{wtype1extpasaasubrets}, $L$ is of type-$2$ $\mathcal{X}_{i}$-extending.

    $\Longleftarrow )$ Let us proceed by induction on the number of classes. It suffices to consider the case $n=2$. Let $a\in L$ and $c\in C(L)$ be such that $a\leq_{e} c$ and $a/0\in \mathcal{X}_{1}\mathcal{X}_{2}$. By the definition of $\mathcal{X}_{1}\mathcal{X}_{2}$, there exists $b\leq a$ such that $a/b \in \mathcal{X}_{2}$ and $b/0 \in \mathcal{X}_{1}$.

Suppose that $\mathcal{X}_{2}$ is closed under initial subintervals. Let $b'$ be a pseudocomplement of $b$ relative to $a/0$. Consider the linear morphism $\pi: a/0 \longrightarrow a/b$ given by $\pi(x)=(x\vee b)\wedge a$. Since $b\wedge b'=0$, we have $\pi(b'/0)=\pi(b')/b=b'/0$. Given that $\mathcal{X}_{2}$ is closed under initial intervals and lattice isomorphisms, we conclude that $b'/0 \in \mathcal{X}_{2}$. Since $b'$ is a pseudocomplement of $b$ in $a/0$, it follows that $b\oplus b^{\prime} \leq_e {a}$, and hence $b\oplus b^{\prime} \leq_{e} c$. By Theorem \ref{sumadirectatipo2},  $L$ is type-$2$ $(\mathcal{X}_{1} \oplus \mathcal{X}_{2})$-extending and $(b'\oplus b)/0 \in \mathcal{X}_{1} \oplus \mathcal{X}_{2}$ and since $a\leq_{e} c$, $b\oplus b'\leq_{e} c$, . Therefore, $c\in D(L)$.

Now, suppose that $\mathcal{X}_{2}$ is closed under quotients. Since $L$ is an idiom, there exists $c^{\prime} \in L$ such that $b\leq_{e} c^{\prime}$ and $c^{\prime}\in C(c/0)$. By Lemma \ref{0.20}, $c^{\prime} \in C(L)$. Since $b/0\in \mathcal{X}_{1}$, we have $c^{\prime}\in D(L)$. Let $c^{\prime\prime}\in L$ be such that $c\oplus c^{\prime\prime}=1$. By modularity, $c=c\wedge (c\oplus c^{\prime\prime})=c\oplus(c \wedge c^{\prime\prime})$ and $c \oplus ((a\vee c)\wedge c^{\prime\prime})= (a\vee c)\wedge (c\oplus c^{\prime\prime})=a\vee c$. Again, using modularity, we have that:
\begin{center}
     $(a\vee c^{\prime})\wedge c^{\prime\prime}/0=((a\vee c^{\prime})\wedge c^{\prime\prime})/((a\vee c^{\prime})\wedge c^{\prime\prime}\wedge c^{\prime})\cong (a\vee c^{\prime})/c^{\prime} \cong a/a\wedge c^{\prime}$.
\end{center}
Given that $a/b\in \mathcal{X}_{2}$, it follows that $a/a\wedge c'\in \mathcal{X}_{2}$, and therefore, $(a\vee c')\wedge c'/0 \in \mathcal{X}_{2}$. Since $a\leq_{e} c$, we have $a\wedge c^{\prime\prime}\leq_{e} c\wedge c'$, so $(a\vee c')\wedge c^{\prime\prime}\leq_{e} c\wedge c^{\prime}$. Since $c=c'\oplus (c\wedge c^{\prime\prime})$, by Lemma \ref{cerrados}, $c\wedge c^{\prime\prime}\in C(L)$. Since $L$ is $\mathcal{X}_{2}$-extending, $c\wedge c^{\prime\prime}\in D(L)$. Furthermore, by modularity, $c\wedge c^{\prime\prime}\in D(c^{\prime\prime}/0)$. Let $t\in L$ be such that $t\oplus (c\wedge c'')=c^{\prime\prime}$, then, by Lemma $0.29$, $c\wedge t=0$ and, by modularity, $c\oplus t=(c'\oplus (c\wedge c''))\oplus t= c'\oplus ((c\wedge c')\oplus t)= c'\oplus c''=1$. Therefore, $L$ is type $\mathcal{X}_{1}\mathcal{X}_{2}$-extending.
\end{proof}
\begin{thm}

   Let $n \in \mathbb{N}$ and $\mathcal{X}_{i}$ ($1 \leq i \leq n$) be a class of lattices such that $\mathcal{X}_{i}$ is closed under initial intervals for each $2 \leq i \leq n$. An idiom $L$ is type-$1$ $(\mathcal{X}_{1}....\mathcal{X}_{n})$-extending if and only if $L$ is type-$1$ $\mathcal{X}_{i}$-extending for all $1 \leq i \leq n$.
\end{thm}
\begin{proof}
    $\Longrightarrow )$ Given that $\mathcal{X}_{i}\subseteq \mathcal{X}_{i}...\mathcal{X}_{n}$ for each $i\leq i \leq n$, according to Lemma \ref{wtype1extpasaasubrets}, $L$ is type-$1$ $\mathcal{X}_{i}$-extending for $1\leq i \leq n$.  
      
      $\Longleftarrow )$ It is enough to prove it for $n=2$. Let $L$ be an idiom such that $L$ is type-$1$ $\mathcal{X}_{i}$-extending with $i=1,2$. Let $a\in L$ be such that $a/0\in \mathcal{X}_{1}\mathcal{X}_{2}$, then there exists $a_{1}\in L$ such that $a_{1}/0 \in \mathcal{X}_{1}$ and $a/a_{1} \in \mathcal{X}_{2}$.  
      Let $b\in P(a)$. We aim to prove that $b\in D(L)$. Since $L$ is an idiom, a pseudocomplement $b'$ relative to $a_1$ exists in $a/0$. Let $\pi: a/0 \longrightarrow a/a_{1}$ be the linear projection morphism given by $\pi(x)=(a_{1}\vee x)$. Taking into account $\pi|_{b^{\prime}/0}$, we find that $k_{\pi_{b^{\prime}/0}} = k_{\pi} \wedge b^{\prime} \leq a_{1} \wedge b^{\prime} = 0$; thus, $\pi(b^{\prime})/a_{1} \cong b^{\prime}/0$.
     Consequently, $b^{\prime}/0 \in \mathcal{X}_{2}$,  because $\mathcal{X}_{2}$ is closed under taking initial intervals and lattice isomorphisms. Given that $b^{\prime}$ is a pseudo-complement of $a_{1}$ in $a/0$, we have $b^{\prime} \oplus a_{1} \leq_{e} a$, which implies $b \in P(b^{\prime} \oplus a_{1})$. According to Theorem \ref{sumadirectatipo1}, $L$ is type-1 $(\mathcal{X}_{1} \oplus \mathcal{X}_{2})$-extending and $a_{1}\oplus b' \in \mathcal{X}_{1} \oplus \mathcal{X}_{2}$, indicating that $b \in D(L)$.
\end{proof}
\begin{cor}\label{corpro}
If $L$ is an idiom,  $\mathcal{S}$ denotes the class of simple lattices, and $\mathcal{L}_{f}$ denotes the class of lattices of finite length, then $L$ is type-$2$ $\mathcal{L}_{f}$-extending if and only if $L$ is type-$2$ $\mathcal{S}$-extending.
 
 Now, $\mathcal{X}_{2}$ is closed under initial intervals, so $\pi(b^{\prime})/a_{1}\in \mathcal{X}_{2}$, and from the previous paragraph, $b^{\prime}/0 \in \mathcal{X}_{2}$. Since $b^{\prime}\oplus a_{1}\leq_{e} a$, we have $b\in P(a_{1}\oplus b^{\prime})$. Finally, since $L$ is of type $(\mathcal{X}_{1}\oplus \mathcal{X}_{2})^{e}$-extending, it follows that $b\in D(L)$.
\end{cor}
\begin{proof}
    $\Longrightarrow )$ This follows from Lemma \ref{wtype1extpasaasubrets} since $\mathcal{S} \subseteq \mathcal{L}_{f}$.
       
$\Longleftarrow )$ Since $\mathcal{S}$ is closed under initial intervals (also, under quotients), by Theorem \ref{sumadirectatipo2}, $L$ is of type $\mathcal{S}^{n}$ for all $n\in \mathbb{N}$. Let $a\in L$ be such that $a/0\in \mathcal{L}_{f}$. Then, there exist $a_{1},...,a_{m}\in L$ such that $0= a_{0} \leq a_{1}\leq a_{2}\leq ... \leq a_{m}=a$, with $a_{i}/a_{i-1}\in \mathcal{S}$ for all $1\leq i \leq m$. Therefore, $a/0\in \mathcal{S}^{n}$, with which $L$ is type-$2$ $\mathcal{L}_{f}$-extending.

\end{proof}
 \begin{cor}
       Let $L$ be a finite modular lattice. $L$ is type-$1$ $\mathcal{S}$-extending if and only if $L$ is an extending lattice.
    \end{cor}
    \begin{proof}
       $\Longleftarrow )$ There's nothing to prove.
       
        $\Longrightarrow )$ Let us observe that $L$ is extending if and only if $L$ is type-1 $\mathcal{L}_{f}$ because any initial interval of $L$ has finite length. By Corollary \ref{corpro}, $L$ is extending.

    \end{proof}
\section{Quasi-continuous lattices}
The following definition naturally generalizes the concept of a quasi-continuous module. Specifically, a left module \(_RM\) is quasi-continuous if and only if its lattice of submodules is also quasi-continuous.

\begin{defn}

A bounded lattice $L$ is quasi-continuous if, for every $a, b \in L$ satisfying $a \wedge b = 0$, there exists $c, d \in L$ such that $c \oplus d = 1$, $a \leq c$, and $b \leq d$.

\end{defn}

\begin{prop}Let $L$ be a bounded modular lattice. Suppose that $L$ satisfies the following conditions:

\begin{enumerate}
\item For every $a \in L$, there exists $d \in D(L)$ such that $a \leq_{e} d$.
\item For all $a,b \in D(L)$ with $a \wedge b = 0$, it follows that $a \oplus b \in D(L)$.\end{enumerate}
Then $L$ is quasi-continuous. 
If $L$ is an idiom, the converse implication holds.
\end{prop}

\begin{proof}
   Let $a, b \in L$ be such that $a \wedge b = 0$. By Condition 1), there exist $c, d \in D(L)$ such that $a \leq_{e} c$ and $b \leq_{e} d$. Then $c \wedge d = 0$, because otherwise, $0\neq c \wedge d \leq c$, $a\leq_e c$ imply $0\neq a\wedge c \wedge d \leq d$.  As $b\leq_e  d$,  then $0\neq b\wedge a\wedge c\wedge d \leq a\wedge b = 0$, a contradiction. By Condition 2, $c \oplus d \in D(L)$. Let $t \in L$ be such that $(c \oplus d) \oplus t = 1$, by Lemma \ref{modular}, $c \oplus (d \oplus t) = 1$ and $a \leq c$, $b \leq d \oplus t$. Therefore, $L$ is quasi-continuous.
     
     Now, assume $L$ is an idiom and quasi-continuous. Let $a \in L$, since $L$ is upper continuous, there exists $b \in P(a)$. Given that $a \wedge b = 0$, there exist $c, d \in L$ such that $c \oplus d = 1$, $a \leq c$, and $b \leq d$.
   Let us demonstrate that $a \leq_{e} c$. Consider $x \leq c$ such that $a \wedge x = 0$. Given that $(a \vee x) \wedge b \leq c \wedge d = 0$, it follows, by Lemma \ref{modular}, that $a \wedge (b \vee x) = 0$. Since $b \in P(a)$, it implies $b \vee x = b$, therefore $x \leq b \wedge c = 0$.
      On the other hand, suppose that $a \wedge b = 0$ and $a, b \in D(L)$. Then there exist $c, d \in D(L)$ such that $a \leq c$, $b \leq d$, and $c \oplus d = 1$. Since $L$ is modular, there exist $a', b' \in L$ such that $a \oplus a' = c$ and $b \oplus b' = d$. Then, by Lemma \ref{modular}, $1 = c \oplus d = (a \oplus a') \oplus (b \oplus b') = (a \oplus b) \oplus (a' \oplus b')$, thus $a \oplus b \in D(L)$.
   \end{proof}
   This definition is presented for modules in \cite{Dogruöz3} page 243.
  
\begin{defn}
   Let \( \mathcal{X} \) be a class of lattices. We say that a bounded lattice \( L \) satisfies \( Q(\mathcal{X}) \), if for each \( a, b \in L \) with \( a/0 \in \mathcal{X} \) and \( a \wedge b = 0 \), there exist \( c, d \in L \) such that \( 1 = c \oplus d \), \( a \leq c \), and \( b \leq d \).
\end{defn}
\begin{rem}
    Let $L$ be a bounded lattice. $L$ is quasi-continuous if and only if $L$ satisfies $Q(\mathcal{L})$.
\end{rem}
\begin{prop}\label{Qxrelativo}
Then,\end{prop}
\begin{proof}
    Let $a\in D(L)$ and $a'\in L$ be such that $a\oplus a'=1$. Consider $c,d\in a/0$ such that  $c/0\in \mathcal{X}$ and  $c\wedge b=0$.
    
Since \( c/0 \in \mathcal{X} \) and \( c \wedge (b \oplus a') = 0 \), by Lemma \ref{modular}, there exist \( d, d' \in L \) such that \( c \leq d \), \( b \oplus a' \leq d' \), and \( d \oplus d' = 1 \).

By modularity, we have \( d' = d' \wedge 1 = d' \wedge (a \oplus a') = a' \oplus (d' \wedge a) \). Therefore, it follows that \( 1 = d \oplus d' = d \oplus (a' \oplus (d' \wedge a)) \). Thus, \( a = a \wedge 1 = a \wedge (d \oplus (a' \oplus (d' \wedge a))) = (a \wedge (d \oplus a')) \oplus (d' \wedge a) \). Here, we conclude that \( c \leq a \wedge (d \oplus a') \) and \( b \leq (d' \wedge a) \). Therefore, we conclude that \( a/0 \) satisfies \( Q(\mathcal{X}) \).

\end{proof}
\begin{lemma}\label{lemmaQx}
   Let $\mathcal{X}$ be a class of lattices and $L$ a bounded modular lattice that satisfies $Q(\mathcal{X})$. Let $a,b \in L$ be such that $a/0 \in \mathcal{X}$ and $b \in P(a)$. Then $1 = a' \oplus b$ with $a \leq_{e} a'$.
\end{lemma}
\begin{proof}
   Since $a\wedge b=0$ and $L$ satisfies $Q(\mathcal{X})$, there exist $a',b' \in L$ such that $a'\oplus b'=1$, $a\leq a'$ and $b\leq b'$. Given that $a\wedge b'=0$ and $b\in P(a)$, we have $b=b'$. By Lemma \ref{independiente}, $(a\oplus b)\leq_{e} (a' \oplus b)=1$. By Lemma \ref{esenmod1}, $(a\oplus b)\wedge a'\leq_{e} 1\wedge a'=a'$. By modularity, $a=(b\wedge a')\vee a=(a\oplus b)\wedge a'\leq_{e} a'.$
\end{proof}
\begin{cor}
    Let $\mathcal{X}$ be a class of lattices and $L$ an idiom. Then for every $a \in L$ with $a/0 \in \mathcal{X}$, there exist $a', b' \in L$ with $a' \oplus b = 1$, with $a \leq_{e} a'$ and $b \in P(a)$.
\end{cor}
\begin{defn}
  Let $\mathcal{X}$ be a class of lattices and $L$ a bounded lattice. We say that $L$ is $\mathcal{X}$-quasi continuous if $L$ satisfies the following properties:
   \begin{enumerate}
\item[$(C1)_{\mathcal{X}}$] For every $a\in L$ with $a/0 \in \mathcal{X}$, there exists $b\in D(L)$ such that $a\leq_{e} b$, and
\item[$(C3)_{\mathcal{X}}$] For all $a,b\in D(L)$ with $a/0 \in \mathcal{X}$ and $a\wedge b=0$, it follows that $a\oplus b\in D(L)$. 
\end{enumerate}

    \end{defn}
    \begin{lemma}\label{C3equiv}
   Let $\mathcal{X}$ be a class of lattices and $L$ a bounded and modular lattice. $L$ satisfies $(C3)_{\mathcal{X}}$ if and only if for all $a, b \in D(L)$ with $a/0 \in \mathcal{X}$ such that $a \wedge b = 0$, there exists $b' \in L$ such that $1 = a \oplus b'$ and $b \leq b'$.
\end{lemma}
\begin{proof} 
$\Longrightarrow )$  Let $a,b \in D(L)$ with $a \in \mathcal{X}$ and $a \wedge b = 0$. Since $L$ satisfies $(C3)_{\mathcal{X}}$, there exists $d \in L$ such that $(a \oplus b) \oplus d = 1$, by Lemma \ref{esenmod1}, $a \oplus (b \oplus d) = 1$, and $b \leq b \oplus d$.

  $\Longleftarrow )$ Let $a,b \in D(L)$ with $a/0 \in \mathcal{X}$ and $a \wedge b = 0$. By hypothesis, there exists $b' \in L$ such that $a \oplus b' = 1$ and $b \leq b'$. Let $c \in L$ be such that $b \oplus c = 1$; then $b' = b' \wedge 1 = b' \wedge (b \oplus c) = b \oplus (c \wedge b')$. Thus, $1 = a \oplus b' = a \oplus (b \oplus (c \wedge b')) = (a \oplus b) \oplus (c \wedge b')$. Therefore, $a \oplus b \in D(L)$. 

\end{proof}
\begin{lemma}\label{quasitype1}
   Let $\mathcal{X}$ be a class of lattices and $L$ an idiom that satisfies $Q(\mathcal{X})$. Then $L$ is type 1 $\mathcal{X}$-extending and $L$ satisfies $(C1)_{\mathcal{X}}$.
\end{lemma}
\begin{proof}

We begin by demonstrating that $L$ is type $1$ $\mathcal{X}$-extending. Let $a \in L$ be such that $a/0 \in \mathcal{X}$ and $b \in P(a)$. Since $L$ satisfies $Q(\mathcal{X})$, there exist $c_{1}, c_{2} \in L$ such that $1 = c_{1} \oplus c_{2}$, $a \leq c_{1}$, and $b \leq c_{2}$. Given $a \wedge c_{2} = 0$, then $b = c_{2}$, and consequently, $b \in D(L)$.

   On the other hand, let $x \in L$ with $x/0 \in L$. Since $L$ is an idiom, there exists $c \in P(x)$. Using an argument similar to the previous paragraph, there exists $d \in D(L)$ with $x \leq d$ and $d \oplus c = 1$. We have $x \leq_{e} d$. To see this, assume $x \wedge t = 0$ for some $t \leq d$. Then, $c \wedge (x \vee t) = 0$. By Lemma \ref{modular}, it follows that $x \wedge (c \vee t) = 0$. Given that $c \in P(x)$, this implies $t = 0$.

\end{proof}
\begin{cor}\label{corQ(X)}
  Let $\mathcal{X}$ be a class of lattices and $L$ an idiom that satisfies $Q(\mathcal{X})$. Then $L$ is $\mathcal{X}$ quasi-continuous.
\end{cor}
\begin{proof}
   This follows from Lemmas \ref{quasitype1} and \ref{C3equiv}.
\end{proof}
\begin{thm} Let X be a class of lattices with $\mathcal{X}=\mathcal{X}^{e}$ and L an idiom. The following are equivalent: 
\begin{enumerate}
  \item L satisfies $Q(\mathcal{X})$. 
  \item L is $\mathcal{X}$-quasi-continuous and type 1 $\mathcal{X}$-extending. 
  \item L is of type 1 and of type 2 $\mathcal{X}$-extending and satisfies $(C3)_{\mathcal{X}}$. 
\end{enumerate} 
\end{thm}
\begin{proof}
 
    \begin{enumerate}\item[] 
      \begin{enumerate}\item[]\item[$1) \Longrightarrow 2)$] It follows from Lemma \ref{quasitype1} and Corollary \ref{corQ(X)}.\item[$2) \Longrightarrow 3)$] Let $a,b \in L$ such that $a/0 \in \mathcal{X}$, $a\leq_{e} b$ and $b\in C(L)$. Given that $\mathcal{X}=\mathcal{X}^{e}$, we have $b\in \mathcal{X}$, and since $L$ satisfies $(C1)_{\mathcal{X}}$, it follows that $b\in D(L)$. Therefore, $L$ is of type $2$ $\mathcal{X}$-extending. Finally, $L$ satisfies $(C3)_{\mathcal{X}}$ for being $\mathcal{X}$-quasi continuous.\item[$3) \Longrightarrow 1)$] Let $a,b \in L$, such that $a/0 \in L$ and $a\wedge b=0$. Since $L$ is an idiom, there exists $c\in P(a)$ such that $b\leq c$. Given that $L$ is of type $1$ $\mathcal{X}$-extending, we have $c\in D(L)$. By Lemma \ref{cerees}, there exists $d\in C(L)$ such that $a\leq_{e} d$, and as $L$ is of type $2$ $\mathcal{X}$- extending, then $d\in D(L)$. As $\mathcal{X}$ is essentially closed, $d/0 \in \mathcal{X}$.
      \end{enumerate} 
      \end{enumerate}
Finally, given that $c \wedge d = 0$, $d/0 \in \mathcal{X}$ and $L$ satisfies $(C3)_{\mathcal{X}}$, there exists $e \in L$ such that $1 = c \oplus (d \oplus e)$. Since $a \leq c$ and $b \leq d \oplus e$, we conclude that $L$ satisfies $Q(\mathcal{X})$.
    \end{proof}

 The following lemma is straightforward.
 
\begin{lemma}\label{Qxconten}
   Let \(\mathcal{X}\) and \(\mathcal{Y}\) be classes of lattices such that \(\mathcal{X} \subseteq \mathcal{Y}\), and let \(L\) be a bounded lattice. If \(L\) satisfies \(Q(\mathcal{Y})\), then \(L\) also satisfies \(Q(\mathcal{X})\).
\end{lemma}
\begin{prop}\label{sumQx} Let $n$ be a positive integer and let $\mathcal{X}_{i}$ be lattice classes for $1\leq i\leq n$. An idiom $L$ satisfies $Q(\mathcal{X}_{1}\oplus \mathcal{X}_{2}\oplus ...\oplus \mathcal{X}_{n})$ if and only if  $L$ satisfies $Q(\mathcal{X}_{i})$ for all $1\leq i\leq n$. 
\end{prop}
\begin{proof}

    $\Longrightarrow )$ It follows from Lemma \ref{Qxconten}.
    
     $\Longleftarrow )$  By induction, it suffices to prove the case $n=2$. Let $a,b \in L$ with $a\in \mathcal{X}_{1}\oplus\mathcal{X}_{2}$ and $a\wedge b=0$. By definition of $\mathcal{X}_{1}\oplus\mathcal{X}_{2}$, there exist $a_{1},a_{2} \in L$ such that $a=a_{1}\oplus a_{2}$ and $a_{1}/0 \in \mathcal{X}_{1}$, $a_{2}/0 \in \mathcal{X}_{2}$.\ Since $L$ is modular, by Lemma \ref{modular}, $a\oplus b=a_{1}\oplus (a_{2}\oplus b)$; as $L$ is $Q(\mathcal{X}_{1})$, there exist $c_{1},c_{2}\in L$ such that $a_{1}\leq c_{1}$, $a_{2}\oplus b\leq c_{2}$ and $c_{1}\oplus c_{2}=1$.\ On the other hand, by Proposition \ref{Qxrelativo}, $c_{2}/0$ satisfies $Q(\mathcal{X}_{2})$, so there exist $d_{1},d_{2}\in c_{2}/0$ such that $c_{2}=d_{1}\oplus d_{2}$, $a_{2}\leq d_{1}$ and $b\leq d_{2}$.\ We conclude that $1=c_{1}\oplus c_{2}=c_{1}\oplus (d_{1}\oplus d_{2})=(c_1{}\oplus d_{1})\oplus d_{2}$ with $a\leq c_{1}\oplus d_{1}$ and $b\leq d_{2}$. Therefore, $L$ satisfies $Q(\mathcal{X}_{1}\oplus\mathcal{X}_{2})$.
\end{proof}
\begin{cor}
   Let $\mathcal{X}$ represent a class of lattices. An idiom $L$ fulfills $Q(\mathcal{X})$ if and only if it fulfills $Q(\mathcal{X}^{\oplus})$.

\end{cor}
\begin{proof}
   This follows from Proposition \ref{sumQx}.
\end{proof}
\begin{thm}
   Let $n \in \mathbb{N}$ and $\mathcal{X}_{i}$ be classes of bounded lattices for $1 \leq i \leq n$. Define $\mathcal{X} = \mathcal{X}_{1} \oplus \ldots \oplus \mathcal{X}_{n}$ and let $L$ be an idiom. Then, $L$ satisfies $Q(\mathcal{X}^{e})$ if and only if $L$ satisfies $Q(\mathcal{X}^{e}_{i})$ for each $1 \leq i \leq n$.
\end{thm}
\begin{proof}
   ($\Longrightarrow$) for each $1 \leq i \leq n$, $\mathcal{X}_{i} \subseteq \mathcal{X}$, which implies $\mathcal{X}_{i}^{e} \subseteq \mathcal{X}^{e}$. By Lemma \ref{Qxconten}, $L$ satisfies $Q(\mathcal{X}_{i}^{e})$ for all $1 \leq i \leq n$.
     $\Longleftarrow )$ We proceed by induction. Clearly, it is enough to prove the result for $n=2$. Let $a, b \in L$ be such that $a/0 \in (\mathcal{X}_{1} \oplus \mathcal{X}_{2})^{e}$ and $a \wedge b = 0$. By Lemma \ref{cerees}, there exists $c \in C(L)$ such that $a \leq_{e} c$; thus $c \wedge b = 0$. Given that $a/0 \in \mathcal{X}_{1} \oplus \mathcal{X}_{2}$, there exist $a_{1}, a_{2} \in L$ such that $a_{1} \oplus a_{2} \leq_{e} a$ and $a_{i}/0 \in \mathcal{X}_{i}$ for $i = 1, 2$.
      By Lemma \ref{cerees}, there exists $a_1' \in C(c/0)$ such that $a_1 \leq_e a_1'$, since $a_1 \leq_e a_1'$ and $a_1 \in \mathcal{X}_1^e$, then $a_1' \in \mathcal{X}_1^e$. By Lemma \ref{cerrados}, $a_1' \in C(L)$. Given that $L$ is an idiom, there exists $b' \in P(a_1')$ such that $a_2 \oplus b \leq b'$, by Lemma \ref{lemmaQx}, there exists $d \in L$ with $a_1' \leq_e d$ and $1 = d \oplus b'$.  Since $a_1' \in C(L)$, we conclude that $d = a_1'$. By modularity, it follows that $c = c \wedge 1 = c \wedge (a_1' \oplus b') = a_1' \oplus (c \wedge b')$ and $(a_1 \oplus a_2) \wedge (c \wedge b') = (a_1 \wedge (c \wedge b')) \oplus a_2 = a_2$ is essential in $c \wedge b'$ by Lemma \ref{esenmod1}. Then $c \wedge b'/0 \in \mathcal{X}_2^e$. Now, by Proposition \ref{Qxrelativo}, $b'/0$ satisfies $Q(\mathcal{X}^e)$; therefore, there exist $e, f \in b'/0$ such that $e \oplus f = b'$, $a \wedge b' \leq e$, and $b \leq f$.
      Finally, note that $1=a_{1}'\oplus b'=a_{1}'\oplus (e\oplus f)=(a_{1}'\oplus e)\oplus f$, where $a\leq c=a_{1}'\oplus (c\wedge b')\leq a_{1}'\oplus e$ and $b\leq f$. Therefore, $L$ satisfies $Q(\mathcal{X}^{e})$.
  
\end{proof}

\noindent\textbf{Jesus Adrian Celis-Gonz\'alez}\\
Departamento de Matemáticas, Facultad de Ciencias,  \\ 
Universidad Nacional Autónoma de México, Mexico City, Mexico. \\ 
\textbf{e-mail:} \textit{celis@ciencias.unam.mx}

\noindent\textbf{Hugo Alberto Rinc\'on-Mej\'ia}\\
Departamento de Matemáticas, Facultad de Ciencias,  \\ 
Universidad Nacional Autónoma de México, Mexico City, Mexico. \\ 
\textbf{e-mail:} \textit{hurincon@gmail.com}\\

\end{document}